\documentclass{amsart}
\usepackage{amsmath,amsthm}
\usepackage{amsfonts,amssymb}
\usepackage{enumerate}
\usepackage[latin1]{inputenc}
\usepackage{accents,color}
\usepackage{graphicx}

\newcommand{\lvt}{\left|\kern-1.35pt\left|\kern-1.3pt\left|}
\newcommand{\rvt}{\right|\kern-1.3pt\right|\kern-1.35pt\right|}

\newtheorem{thm}{Theorem}
\newtheorem{cor}[thm]{Corollary}

\newtheorem{prob}[thm]{Problem}

\newtheorem{THEO}{Theorem}

\theoremstyle{remark}

\graphicspath{{./}}

\begin{document}

\title[Markov's theorem]{Monotonicity of zeros of polynomials orthogonal with respect to a discrete measure}

\author{Dimitar K. Dimitrov}
\address{Departamento de Matem\'atica Aplicada\\
 IBILCE, Universidade Estadual Paulista\\
 15054-000 Sa\~{o} Jos\'e do Rio Preto, SP, Brazil.}
 \email{dimitrov@ibilce.unesp.br}

\keywords{orthogonal polynomials, zeros, monotonicity, discrete measure}
\subjclass[2010]{33C45, 26C10}

\begin{abstract}
We prove that all zeros of the polynomials orthogonal with respect to a measure
$d \mu(x;a) = d \mu(x) + M \delta(x-a)$, where $d\mu$ is a nonatomic positive Borel 
measure and $M>0$, are increasing functions of the mass point $a$. Thus we solve partially an 
open problem posed by Mourad Ismail. 
\end{abstract}

\dedicatory{Dedicated to Mourad Ismail who inspired my interest in zeros of 
orthogonal polynomials.} 

\maketitle

\section{Introduction}
\setcounter{equation}{0}

Let $d\alpha(x;\tau)$ be a family of positive Borel measures depending on the parameter $\tau \in (c,d)$ and $\{ p_n(x;\tau)\}$ be the sequence 
of polynomials orthogonal with respect to $d\alpha(x;\tau)$.  For a given $n\in \mathbb{N}$, denote by $x_{nk}(\tau)$, $k=1,\ldots,n$, the zeros 
of $p_n(x;\tau)$. Markov \cite{Markov} and Stieljes \cite{Sti} were the first who studied $x_{nk}(\tau)$ as functions of the the parameter $\tau$ 
in two fundamental papers, both published in 1886. Markov's beautiful result \cite{Markov} (see also \cite[Theorem 7.1.1]{Ismail} and \cite[Theorem 6.12.1]{Szego})  
concerns the case when the measure is absolutely continuous on $(a,b)$, that is $d\alpha(x;\tau) = \omega(x;\tau) dx$, where,  for any $\tau\in (c,d)$,  $\omega(x;\tau)$ is a weight 
function supported on $(a,b)$.

\begin{THEO} 
\label{A}
Let $\{ p_n(x;\tau) \}$ be orthogonal with respect to $d \alpha(x;\tau)$,
$$
d \alpha(x;\tau) = \omega(x;\tau) d \alpha(x),
$$
on an interval $I=(a,b)$ and assume that $\omega(x;\tau)$ is positive and has continuous 
first derivative with respect to $\tau$ for $x\in (a,b)$, $\tau \in T=(\tau_1,\tau_2)$. Furthermore, 
assume that 
$$
\int_{a}^{b} x^j \omega_\tau(x;\tau)\, d\alpha(x),\ \ j=0, 1, \ldots, 2n-1,
$$
converge uniformly for $\tau$ in every compact subinterval of $T$. Then the zeros of $p_n(x;\tau)$
are increasing (decreasing) functions of $\tau$, $\tau \in T$, if $\partial \{ \ln \omega(x;\tau)\} /\partial \tau$ 
is an increasing (decreasing) function of $x$, $x\in I$.
\end{THEO}

The lack of similar results in the case when the measure contains discrete masses motivated Ismail to formulate the following problem (see \cite[Problem[24.9.1]{Ismail}):
\begin{prob} 
Extend Theorem A to the case when 
$$ 
d \alpha(x;\tau) = \omega(x;\tau) d x + d \beta(x;\tau),
$$
where $\beta(x;\tau)$ is a jump function or a step function.
\end{prob}

Then Ismail emphasises the importance of a result when only the discrete part of the measure depends on the parameter: 
``The case of purely discrete measures is of particular interest so we pose the problem
of finding sufficient conditions on $d \beta(x;\tau)$ to guarantee the monotonicity of the
zeros of the corresponding orthogonal polynomials when the mass points depend on
the parameter $\tau$.''

The above problem, as stated, is challenging because of its generality. However, the absolutely continuous part and the 
discrete part of the measure may force the zeros to move either in the same or in the opposite direction. In the latter situation 
the influence of each part of of the measure may be rather complex so that one could hardly expect reasonable 
sufficient conditions. The above comment of Ismail might have been motivated by the latter argument as well as 
by the fact that the influence of the absolutely continuous part is described in Theorem \ref{A}. Therefore, we concentrate on 
the discrete part only. In other words, from now on we consider the case
\begin{equation}
d \alpha(x;\tau) = \omega(x)\, dx + d \beta(x;\tau).
\label{dal}
\end{equation}
In general, 
\begin{equation}
 d \beta(x;\tau) = \sum_{j=0}^{\kappa} M_j(\tau)\, \delta(x-a_j(\tau)).
 \label{dbeta} 
\end{equation}
where $M_j(\tau)>0$, $a_j(\tau) \in \mathbb{R}$ and $\delta$ stands for the Dirac delta. Each term depends on both  $M_j(\tau)$ 
and $a_j(\tau)$. When all positions $a_j(\tau)$ are fixed, that is $a_j(\tau) = a_j$, and only $M_j(\tau)$ depend on $\tau$, one may 
simply apply Theorem A with an obvious modification: extend $M_j(\tau)$ to a positive function $M_j(x;\tau)$ on $(a,b)$ and look at 
$\alpha(x)$ as a Stietjes distribution with jumps at $a_j$. The corresponding results about the monotonic behaviour of the zeros of 
the classical discrete orthogonal polynomials can be found in \cite[Theorem A and Theorem 2.1]{MCOMP}.

Observe also that the influence of the distinct discrete parts in (\ref{dbeta}) could be pretty complex, so that we concentrate on a discrete 
measure with a single mass, that is we reduce the problem to the case
 $$
 d \beta(x;\tau) = \sum_{j=0}^{\kappa} M_j(\tau)\, \delta(x-a_j(\tau)) 
 $$

Therefore, we arrive at the question about the monotonicity of zeros of polynomials orthogonal with respect to the location of a single mass point. Thus we 
let 
\begin{equation}
\label{dmat}
d \mu(x;\tau) = d \mu(x) + M \delta(x-a(\tau)),
\end{equation}
where $d\mu$ is a nonatomic positive Borel measure,  and prove:
\begin{thm} 
If $d \mu(x;\tau)$ is defined by (\ref{dmat}) and $a(\tau)$ is an increasing function of $\tau$, then all zeros $x_{nk}(\tau)$ of the polynomials 
$p_n(x;\tau)$, orthogonal with respect to $d \mu(x;\tau)$, are increasing functions of $\tau$.     
\end{thm}

A simpler equivalent statement is 

\begin{cor} 
\label{corr}
If $d \mu(x;a) = d \mu(x) + M \delta(x-a)$, then all zeros $x_{nk}(a)$ of the polynomials $p_n(x;a)$, orthogonal with respect to $d \mu(x;a)$,
are increasing functions of $\tau$.   
\end{cor}

Succinctly, our main result states that a moving mass point of a measure ``pulls'' all the zeros of the corresponding orthogonal polynomials towards the same direction it moves.

\section{Proof} 

We prove Corollary \ref{corr}.

The idea of the proof is rather simple and natural. One needs to approximate Dirac's delta by the normal distribution centred at $a$, consider $a$ as a parameter 
and apply the classical Markov's theorem, that is, Theorem A to the normal distribution. Thus, the proof is a formalisation of this straightforward idea. 

First we prove the statement in the case when $d \mu(x)$ is an absolutely continuous measure, $d \mu(x)=\omega(x) dx$. 

Let 
$$
\mathcal{N}(x;a,\gamma) = \frac{1}{\sqrt{\pi} \gamma} e^{-((x-a)/\gamma)^2}.
$$
Then 
$$
\int_{-\infty}^{\infty} \mathcal{N}(x;a,\gamma)\, dx = 1\ \ \mathrm{for\ every}\ a\in \mathbb{R}\ \ \mathrm{and}\ \ \gamma>0.
$$
It is well known that $\mathcal{N}(x;a,\gamma)$ converges to $\delta(x-a)$ in the sense of distributions when $\gamma$ goes to zero with positive values. Let $p_n(x;a,\gamma)$ 
be the polynomials orthogonal with respect to $d \mu(x;a,\gamma) = \left\{ \omega(x) + M\, \mathcal{N}(x;a,\gamma) \right\} dx$. Rather straightforward calculations show that for every $a\in \mathbb{R}$ and 
$\gamma>0$ the moments $m_k(a,\gamma)$ of $\mathcal{N}(x;a,\gamma)  dx$ are explicitly given by 
\begin{eqnarray*}
m_k(a,\gamma) & = & \int_{-\infty}^{\infty} x^k\, \mathcal{N}(x;a,\gamma)\, dx \\
\ & = & \frac{1}{\sqrt{\pi}} \sum_{j=0}^{k} {k \choose j} \frac{1+(-1)^j}{2}\, \Gamma \left( \frac{j+1}{2} \right)\, a^{k-j}\, \gamma^j.
\end{eqnarray*}
Then, for any fixed $a\in \mathbb{R}$ and $k\in \mathbb{N}$, $m_k(a,\gamma) \rightarrow a^k$ as $\gamma \rightarrow +0$. On the other hand, the $k$-th moment of $\delta(x-a)$ is exactly $a^k$. 
Therefore the moments of $d \mu(x;a,\gamma)$ converge to the moments of $d \mu(x;a)$ as $\gamma \rightarrow + 0$. It follows from the representation of the orthogonal polynomials 
as quotients of determinants involving the moments of the measure  (see \cite[(2.26)]{Szego} or \cite[(2.1.6)]{Ismail}) that, for any $n\in \mathbb{N}$ and $a \in \mathbb{R}$, the Taylor coefficients of the 
polynomial $p_n(x;a,\gamma)$ converge to those of  $p_n(x;a)$ when $\gamma$ converges to zero with positive values. This immediately yields that $p_n(x;a,\gamma)$ converges locally uniformly 
to $p_n(x;a)$.  Then a well-known theorem of Hurwitz (see \cite[Theorem 14.3.2]{Hille}) implies  that the zeros $x_{nk}(a,\gamma)$ of $p_n(x;a,\gamma)$ converge to the zeros $x_{nk}(a)$ of $p_n(x;a)$ when $\gamma \rightarrow +0$. 
Moreover, by the implicit function theorem $x_{nk}(a,\gamma)$  are smooth functions of $\gamma$, for $\gamma\in (0,\infty)$.

Let us fix $\gamma>0$ and consider the behaviour  of  $x_{nk}(a,\gamma)$ as functions of $a$.  It is clear that 
\begin{eqnarray*}
\mathcal{N}_a(x;a,\gamma) & = & \frac{\partial}{\partial a} \mathcal{N}(x;a,\gamma) \\
\ & = & \frac{2}{\sqrt{\pi} \gamma^3}\ (x-a)\ e^{-((x-a)/\gamma)^2}
\end{eqnarray*}
is a continuous function of $a \in \mathbb{R}$ as well as that the moments of $\mathcal{N}_a(x;a,\gamma) $ are integrals that converge uniformly for $a$ in any subinterval of the real line.  
Since 
$$
\frac{\partial}{\partial a}  \ln \mathcal{N}(x;a,\gamma) = \frac{2 (x-a)}{\gamma^2}
$$
is an increasing function of $x$, then, by Theorem A, all $x_{nk}(a,\gamma)$  are increasing functions of $a$. 
Taking a limit, as $\gamma \rightarrow +0$, we conclude that $x_{nk}(a)$ are increasing functions of $a$ too.
This completes the first proof in the case when $d \mu(x)$ is absolutely continuous. 

Since a nonatomic Borel measure measure contains an absolutely continuous part and an enumerable quantity of jumps, we may approximate 
the corresponding delta functions either by normal distributions or even by $C^\infty$ compactly supported test functions (see \cite{Gel}) and then we proceed as 
above.  

The main result may be applied to deduce the monotonicity of zeros of sequences of polynomials orthogonal with respect to a mesure of the form (\ref{dal}), 
where $\omega$ is associated with the classical sequences of orthogonal polynomials (see \cite{CMC, CMVQ}).

\end{document}